\documentclass[a4paper]{article}
\usepackage{amsmath,amssymb,fancybox,longtable,graphics,amscd,multirow}
\usepackage[dvipdfm]{graphicx}
\usepackage{caption,setspace}
\newtheorem{thm}{\textsc{Theorem}}[section]

\newtheorem{defn}{\textsc{Definition}}[section]

\def\QED{$\Box$}

\def\mbi#1{\boldsymbol{#1}} 

\def\Pic{\mathop{\mathrm{Pic}}\nolimits}

\def\rk{\mathop{\mathrm{rank}}\nolimits}

\begin{document}
\title{Polytope duality for families of $K3$ surfaces associated to transpose duality}
\author{Makiko Mase}
\date{\small \begin{tabular}{l} Key Words: $K3$ surfaces, toric varieties \end{tabular}}
\maketitle
\begin{abstract}
We consider whether or not transpose-dual pairs, which is a Berglund--H$\ddot{\textnormal{u}}$bsch mirror studied by Ebeling and Ploog \cite{EbelingPloog},  extend to a polytope duality that has a potential to be lattice dual. 
\end{abstract}
\section{Introduction} \label{Introduction}
Isolated singularities in $\mathbb{C}^3$ are classified by Arnold \cite{Arnold75} among which there are classes called {\it bi}modal and {\it uni}modal. 
Our notation follows that of Arnold's. 
Not only the classification, Arnold also finds that there is a duality among unimodal singularities that is called {\it Arnold's strange duality}. 
The duality is also related to toric geometry and lattice theory. 
Ebeling and Ploog \cite{EbelingPloog} find an analogous duality concerning bimodal and other singularities, which is actually a Berglund--H$\ddot{\textnormal{u}}$bsch mirror. 

Batyrev's proposal \cite{BatyrevMirror} of polar duality of {\it reflexive} polytopes gives a breakthrough in a construction of mirror partner for {\it toric} Calabi-Yau hypersurfaces and later complete intersections. 

Being origined in physics, there appear a numerical meanings of ``mirror'' such as cohomological mirror, among which in this article we focus on a relation between Ebeling and Ploog's transpose duality and Batyrev's polytope duality associating with bimodal singularities in some manner. 

In a series of recent studies, it is concluded that transpose-dual pairs 
$(Q_{12},\, E_{18})$, $(Z_{1,0},\,  E_{19})$, $(E_{20},\, E_{20})$, $(Q_{2,0},\, Z_{17})$, $(E_{25},\, Z_{19})$, $(Q_{18},\, E_{30})$
of singularities can extend to a lattice duality by the author \cite{Mase15} following an extension to polytope duality by the author and Ueda \cite{MU}. 
However, those pairs in the list $(*)$ below fail to extend to a lattice duality in spite of the fact that they are polytope dual.  
\[
(*)
\begin{array}{l}
(Z_{13},\, J_{3,0}),\, 
(Z_{1,0},\, Z_{1,0}),\, 
(Z_{17},\, Q_{2,0}),\,  
(U_{1,0},\, U_{1,0}),\, 
(U_{16},\, U_{16}), \\
(Q_{17},\, Z_{2,0}),\,  
(W_{1,0},\, W_{1,0}),\, 
(W_{17},\, S_{1,0}),\,  
(W_{18},\, W_{18}),\,  
(S_{17},\, X_{2,0}). 
\end{array}
\]
More precisely, for each pair one obtains in \cite{MU} reflexive polytopes $\Delta_{[MU]}$ and $\Delta_{[MU]}'$ satisfying that the polar dual of $\Delta_{[MU]}$ is isomorphic to $\Delta_{[MU]}'$ and that $\Delta_{[MU]}$ and $\Delta_{[MU]}'$ respectively contains the Newton polytope of a compactification polynomial of the defining polynomial of singularities. 
Despite this fact it is concluded in \cite{Mase15} that the corresponding pairs of families $\mathcal{F}_{\Delta_{[MU]}}$ and $\mathcal{F}_{\Delta_{[MU]}'}$ of $K3$ surfaces are not lattice dual, that is, the Picard lattices $\Pic(\Delta_{[MU]})$ and $\Pic(\Delta_{[MU]}')$ of these families do not satisfy an isometry $\Pic(\Delta_{[MU]})_{\Lambda_{K3}}^\perp\simeq U\oplus \Pic(\Delta_{[MU]}')$. 
Moreover, for these pairs we can observe that the restriction map $H^{1,1}(\widetilde{\mathbb{P}_{\Delta_{[MU]}}},\,\mathbb{Z}) \to H^{1,1}(\widetilde{Z},\,\mathbb{Z})$ for the minimal model of any generic member $Z\in\mathcal{F}_{\Delta_{[MU]}}$ is not surjective. 

The aim of the study is to consider the following problem arisen by Professor Ashikaga's question: \\

\noindent
\begin{problem}\, 
Let $((B,\, f),\, (B',\, f'))$ be a transpose-dual pair in the list $(*)$ together with their defining polynomials $f$ and $f'$. 
Determine whether or not it is possible to take polynomials $F$ and $F'$ that are respectively compactifications of $f$ and $f'$, and a reflexive polytope $\Delta$ such that the following condition $(**)$ holds: 
\[
(**)\qquad \Delta_F\subset \Delta,\,  \Delta_{F'}\subset \Delta^*, \quad \textrm{and} \quad L_0(\Delta)=0. 
\]
Here, $\Delta_F$ and $\Delta_{F'}$ denote respectively the Newton polytopes of $F$ and of $F'$, and $\Delta^*$ is the polar dual polytope of $\Delta$. 
\end{problem} \\

The main theorem of this paper is stated as follows: \\

\noindent
{\bf Main Theorem.} (Theorem \ref{MainThmPf})\, 
{\it For each of the following pairs 
\[
(Z_{1,0},\, Z_{1,0}),\, 
(U_{1,0},\, U_{1,0}),\, 
(Q_{17},\, Z_{2,0}),\,  
(W_{1,0},\, W_{1,0}), 
\]
there exist compactifications $F,\, F'$ and reflexive polytopes $\Delta$ and $\Delta'$ such that 
\[
(**)\qquad \Delta^* \simeq \Delta',\, \Delta_F\subset \Delta,\,  \Delta_{F'}\subset \Delta', \quad \textit{and} \quad \rk L_0(\Delta)=0
\]
hold. 
Moreover, $\rho(\Delta)+\rho(\Delta')=20$. 
}\\

It can be conjectured that there do not exist reflexive polytopes for pairs $(Z_{13},\, J_{3,0})$, $(Z_{17},\, Q_{2,0})$, $(U_{16},\, U_{16})$, $(W_{17},\, S_{1,0})$, $(W_{18},\, W_{18})$, $(S_{17},\, X_{2,0})$ of singularities satisfying the condition $(**)$. 
We leave the judgement about this conjecture to a further study in the furure.  

Section \ref{Preliminary} is devoted to recall some facts as to a polytope duality associated to singularities. 
The proof of the main theorem is given in Theorem \ref{MainThmPf} in section \ref{MainThm}, where we explicitely give compactifications and reflexive polytopes for these pairs. \\

\noindent
\begin{ackn}\\
{\rm {
The author would be grateful to Professor T. Ashikaga for his question of the problem that motivated this article after publication of \cite{Mase15}, to Professor N. Aoki who was  reading through the first draft carefully and making many helpful suggestions, and to Professor M. Kobayashi for his comments and encouragement. 
}}
\end{ackn}
\section{Preliminary} \label{Preliminary}
Recall that a {\it Gorenstein $K3$ surface} is a compact complex connected $2$-dimensional algebraic variety $S$ with at most $ADE$ singularities satisfying $K_S\sim 0$ and $H^1(S,\, \mathcal{O}_S)=0$. 
If a Gorenstein $K3$ surface is nonsingular, we simply call it a {\it $K3$ surface}. 

Let $M\simeq\mathbb{Z}^3$ be a $3$-dimensional lattice and $N=\mathrm{Hom}_{\mathbb{Z}}(M,\, \mathbb{Z})\simeq\mathbb{Z}^3$ the dual of $M$ with a natural pairing $\langle \, , \, \rangle : N\times M \to \mathbb{Z}$. 
Let $\Delta$ be a $3$-dimensional polytope, that is, $\Delta$ is a convex hull of finitely-many points in $M\otimes_{\mathbb{Z}}\mathbb{R}$. 
The associated toric $3$-fold is denoted by $\mathbb{P}_\Delta$. 
The {\it polar dual} $\Delta^*$ of $\Delta$ is defined by
\[
\Delta^* = \left\{ y\in N\otimes_{\mathbb{Z}}\mathbb{R} \, |\, \langle y,\, x\rangle \geq -1 \quad \textnormal{for all} \quad x\in\Delta\right\}. 
\]

Let us recall a toric description of weighted projective spaces. 
Let $\mbi{a}=(a_0, a_1, a_2, a_3)$ be a well-posed quadruple of natural numbers and $d=a_0+a_1+a_2+a_3$. 
Define a $3$-dimensional lattice $\tilde{M}$ by 
\[
\tilde{M} := \left\{ (i, j, k, l)\in\mathbb{Z}^4\, | \, a_0i+ a_1j+ a_2k+ a_3l \equiv 0\mod d\right\} \simeq \mathbb{Z}^3. 
\]
Note that the lattice $\tilde{M}$ is one-to-one corresponding to the set of monomials of weighted degree $d$: indeed, for each $(i, j, k, l)\in\tilde{M}$, a monomial $X_0^iX_1^jX_2^kX_3^l$ is of weighted degree $d$. 
Here, the weight of $X_i$ is $a_i$ for $i=0,1,2,3$. 
Besides, by letting $\Delta_{\mbi{a}}$ be a convex hull of all primitive lattice vectors in $\tilde{M}$, the associated projective toric $3$-fold is the weighted projective space of weight $\mbi{a}$. 

The introduction of {\it reflexive polytope} in \cite{BatyrevMirror} is motivated by mirror symmetry. 
\begin{defn}\textnormal{(\cite{BatyrevMirror})}
Let $\Delta$ be an integral polytope that contains the origin in its interior. 
The polytope $\Delta$ is called {\rm reflexive} if its polar dual $\Delta^*$ is also integral. 
\end{defn}

Not only in a context of mirror, this notion is basically friendly with $K3$ surfaces as follows: 
\begin{thm}\textnormal{(\cite{BatyrevMirror})} \label{ReflexiveThm}
Let $\Delta$ be a $3$-dimensional polytope. \\
$(1)$\, The followings are equivalent: 
\begin{itemize}
\item[$(i)$] The polytope $\Delta$ is reflexive. 
\item[$(ii)$] The toric $3$-fold $\mathbb{P}_\Delta$ is Fano, in particular, general anticanonical members of $\mathbb{P}_\Delta$ are Gorenstein $K3$. 
\end{itemize}
$(2)$\, General anticanonical members of $\mathbb{P}_\Delta$ are simultaneously resolved by a toric (crepant) desingularization of $\mathbb{P}_\Delta$ to be $K3$ surfaces. 
\end{thm}

Denote for a reflexive polytope $\Delta$ by $\mathcal{F}_\Delta$ a family of (Gorenstein) $K3$ surfaces parametrised by the complete anticanonical linear system $|{-}K_{\mathbb{P}_\Delta}|$. 
For a member $Z$ in $\mathcal{F}_\Delta$, denote by $\tilde{Z}$ and $\widetilde{\mathbb{P}_{\Delta}}$ the minimal models in a cause of the simultaneous resolution. 

In the article, we define that a member $Z\in\mathcal{F}_\Delta$ is {\it generic} if the following two conditions are satisfied: 
\begin{itemize}
\item[(1)] $Z$ is $\Delta$-regular. (See \cite{BatyrevMirror} for detail)
\item[(2)] The Picard group of $\widetilde{Z}$ is generated by irreducible components of the restrictions of the generator of the Picard group of $\widetilde{\mathbb{P}_{\Delta}}$. 
\end{itemize}

It is proved in \cite{BatyrevMirror} that $\Delta$-regularity is a general condition. 
The condition (2) is also a general condition. 
Note that all Picard lattices of the minimal models of any generic members are isometric.
\begin{defn}
$(1)$\, The {\rm Picard lattice $\Pic(\Delta)$ of the family $\mathcal{F}_\Delta$} is the Picard lattice of the minimal model of a generic member. \\
$(2)$\, $\rho(\Delta):=\rk\Pic(\Delta)$ is called the {\rm Picard number of the family $\mathcal{F}_\Delta$}. \\
$(3)$\, Let $r:H^{1,1}(\widetilde{\mathbb{P}_{\Delta}},\, \mathbb{Z}) \to H^{1,1}(\tilde{Z},\, \mathbb{Z})$ be the restriction mapping of the cohomology group. 
The cokernel of $r$ is denoted by $L_0(\Delta)$. 
\end{defn}

In \cite{MU}, a notion of transpose duality \cite{EbelingPloog} for singularities is extended to a {\it polytope duality} in the sense of the following theorem : 
\begin{thm}{\rm (\cite{MU})}
Let $((B,\, f),\, (B',\, f'))$ be a transpose-dual pair together with their defining polynomials $f$ and $f'$ that are respectively compactified to polynomials $F$ and $F'$. 
Then, there exist reflexive polytopes $\Delta_{[MU]}$ and $\Delta'_{[MU]}$ such that 
\[
\Delta_{[MU]}^*\simeq \Delta'_{[MU]},\quad \Delta_F\subset \Delta_{[MU]}, \quad \textit{and} \quad \Delta_{F'}\subset \Delta'_{[MU]}. 
\]
\end{thm}

However, it depends on the pairs that whether or not $\rk L_0(\Delta_{[MU]})=0$ holds. 
In section \ref{MainThm}, we shall show that some pairs in the list $(*)$ do have this property. 

We end this section by giving formulas that are needed in the proof of the main theorem. 
See \cite{Kobayashi} for details. 
For a $3$-dimensional reflexive polytope $\Delta$, denote by $\Delta^{[1]}$ the set of all edges of $\Delta$, and for an edge $\Gamma\in\Delta^{[1]}$, the dual edge in the polar dual polytope $\Delta^*$ is denoted by $\Gamma^*$. 
The number of lattice points on an edge $\Gamma$ is denoted by $l(\Gamma)$, whilst $l(\Gamma)-2$ by $l^*(\Gamma)$. 
We have 
\begin{eqnarray}
\rk L_0(\Delta) & = & \sum_{\Gamma\in\Delta^{[1]}} l^*(\Gamma) l^*(\Gamma^*). \label{SurjectionCriterion}\\
\rho(\Delta) & = & \sum_{\Gamma\in\Delta^{[1]}} l(\Gamma^*)-3. \label{Rho}
\end{eqnarray}
Note that $\rk L_0(\Delta)=\rk L_0(\Delta^*)$ by the formula. 

\section{Main result} \label{MainThm}
The chief aim of this section is to prove the following statements. 
\begin{thm}\label{MainThmPf}
For pairs $(B,\, B')$ of singularities, if one takes compactifications $F,\, F'$ as in Table \ref{ExcPairsPolynom}, and polytopes $\Delta,\, \Delta'$ as in Table \ref{ExcPairsPolytope}, then, 
\begin{itemize}
\item[$(i)$] $\Delta$ and $\Delta'$ are reflexive, 
\item[$(ii)$] $\Delta^*$ is isomorphic to $\Delta'$ up to lattice isometry of $\mathbb{Z}^3$, 
\item[$(iii)$] $\Delta_F\subset \Delta,$ and  $\Delta_{F'}\subset \Delta'$ hold, and
\item[$(iv)$] $\rk L_0(\Delta)=0$. 
\end{itemize}
Moreover, $\rho(\Delta)+\rho(\Delta') = 20$. 
\[
\begin{array}{cccc} 
B & F & F' & B' \\ 
\hline
\hline
Z_{1,0} & X^5Y+XY^3+Z^2+W^{10}X^2 &  X^5Y+XY^3+Z^2+W^{14}  & Z_{1,0}\\
\hline
U_{1,0} & X^3Y+Y^2Z+Z^3+WX^4 & XZ^3+X^2Y+Y^3+W^9  & U_{1,0}  \\
\hline
Z_{2,0} & X^5Z+XY^3+Z^2+W^7Y & X^5Y+WY^3+XZ^2+W^7  & Q_{17}\\
\hline
W_{1,0} & X^6+Y^2Z+Z^2+W^6Z & X^6+Y^2Z+Z^2+W^{12}  & W_{1,0}\\
\hline
\end{array}
\]
\begingroup
\captionof{table}{Compactifications of singularities}\label{ExcPairsPolynom}
\endgroup
\[
\begin{array}{cccccc} 
B & \textnormal{vertices of }\Delta & \textnormal{vertices of }\Delta' & B' \\ 
\hline
\hline
Z_{1,0} & \left\{ \begin{array}{l} (-1,0,1), (-1,0,0), \\ (0,1,-1), (2,3,-1), \\ (2,2,-1), (1,-1,-1), \\ (0,-1,-1)\end{array}\right\} & \left\{ \begin{array}{l}(0,2,-1), (-1,1,-1), \\ (-1,-1,-1), (5,-1,-1), \\ (4,0,-1), (1,0,0), \\ (-1,-1,1) \end{array}\right\} & Z_{1,0}\\
\hline
U_{1,0} & \left\{ \begin{array}{l}(-1,0,2), (0,1,0), \\  (1,2,-1), (1,1,-1), \\ (0,-1,0), (0,-1,-1)\end{array}\right\} & \left\{ \begin{array}{l}(1,0,-1), (0,-1,-1), \\ (-1,-1,-1), (-1,2,-1), \\ (1,2,-1), (1,0,1), \\ (0,-1,2), (-1,-1,2) \end{array}\right\} & U_{1,0}  \\
\hline
Z_{2,0} & \left\{ \begin{array}{l}(-1,-1,2), (0,-1,0), \\  (1,-1,0), (1,-1,1), \\ (1,2,-3), (0,0,-1) \end{array}\right\}& \left\{ \begin{array}{l}(-1,2,-1), (-1,-1,1), \\ (-1,-1,-1), (6,-1,-1), \\ (2,1,-1), (0,-1,1) \end{array}\right\}& Q_{17}\\
\hline
W_{1,0} & \left\{ \begin{array}{l}(-1,0,1), (-1,0,0),  \\ (1,2,-1), (2,3,-1), \\ (0,-1,0) \end{array}\right\} & \left\{ \begin{array}{l}(-1,-1,-1), (5,-1,-1), \\ (1,3,-1), (-1,3,-1), \\ (-1,-1,1) \end{array}\right\} & W_{1,0}\\
\hline
\end{array}
\]
\begingroup
\captionof{table}{Polytopes that make the pairs polytope dual}\label{ExcPairsPolytope}
\endgroup
\end{thm}
{\sc Proof.} 

{\bf $Z_{1,0}$ case.}
The defining polynomials of singularities $B=Z_{1,0}$ and $B'=Z_{1,0}$ are the same $f=f'=x^5y +xy^3 + z^2$. 

Take a compactification of $f$ as $F=W^{10}X^2 + X^5Y +XY^3 + Z^2$ in the weighted projective space $\mathbb{P}(1,2,4,7)$. 
Note that $F$ is a different compactification from the one in \cite{EbelingPloog}. 

Take a compactification of $f'$ as $F'=W^{14} + X^5Y +XY^3 + Z^2$ in the weighted projective space $\mathbb{P}(1,2,4,7)$. 
Note that $F'$ is the same compactification as in \cite{EbelingPloog}. 

The polytope $\Delta$ contains the Newton polytope of $F$: indeed, by taking a basis $\mbi{e}_1=(-6,1,1,0)$,\, $\mbi{e}_2=(2,1,-1,0)$,\, $\mbi{e}_3=(-7,0,0,1)$ for $\mathbb{R}^3$, one can see that monomials $W^{10}X^2,\, X^5Y,\, XY^3,\, Z^2$ are respectively corresponding to vertices 
\[
(0,1,-1),\, (2,2,-1),\, (1,-1,-1),\, (-1,0,1). 
\]

The polytope $\Delta'$ contains the Newton polytope of $F'$: indeed, by taking a standard basis $\mbi{e}_1'=(-2,1,0,0)$,\, $\mbi{e}_2'=(-4,0,1,0)$,\, $\mbi{e}_3'=(-7,0,0,1)$ for $\mathbb{R}^3$, one can see that monomials $W^{14}$,\, $X^5Y$,\, $XY^3$,\, $Z^2$ are respectively corresponding to vertices 
\[
(-1,-1,-1),\, (4,0,-1),\, (0,2,-1),\, (-1,-1,1). 
\]

The dual polytope $\Delta'^*$ of $\Delta'$ is a convex hull of vertices 
\[
(0,0,1),\, (-1,-2,-3),\, (-1,-3,-5),\, (1,-1,-1),\, (1,0,0),\, (0,1,0),\, (-1,-1,-3)
\]
that is mapped to isomorphically from $\Delta$ by a transformation of $\mathbb{R}^3$ by the matrix 
\[
M:=\left(
\begin{array}{ccc}
 1 & 2 & 3 \\
 0 & -1 & -1 \\
 1 & 2 & 4 \\
\end{array}
\right)
\]
that is, $(x, y, z)M=(x', y', z')$ for $(x,y,z)\in\Delta$ and $(x',y',z')\in\Delta'$. 

Therefore, $\Delta$ and $\Delta'$ are reflexive and the pair is polytope dual. 

By the formula (\ref{SurjectionCriterion}), one gets $\rk L_0(\Delta)=\rk L_0(\Delta^*)=0$ because for all edges in $\Delta$ satisfy $l^*(\Gamma)l^*(\Gamma^*)=0$. 
In fact, at least either $\Gamma$ or $\Gamma^*$ has no lattice points in its interior. 

By the formula (\ref{Rho}), one can compute that 
\[
\rho(\Delta) = 17-3 = 14, \quad \rho(\Delta^*) = 9-3 = 6
\]
thus one has
\[
\rho(\Delta) +\rho(\Delta^*) = 20. 
\]
\\

{\bf $U_{1,0}$ case.} 
The defining polynomials of singularities $B=U_{1,0}$ and $B'=U_{1,0}$ are $f=x^3y +y^2z + z^3,\, f'=x'z'^3 +x'^2y' + y'^3$, respectively.  

Take a compactification of $f$ as $F=WX^4 + X^3Y +Y^2Z + Z^3$ in the weighted projective space $\mathbb{P}(1,2,3,3)$. 
Note that $F$ is a different compactification from the one in \cite{EbelingPloog}. 

Take a compactification of $f'$ as $F'=W'^9 + X'Z'^3 +X'^2Y' + Y'^3$ in the weighted projective space $\mathbb{P}(1,3,3,2)$. 
Note that $F'$ is the same compactification as in \cite{EbelingPloog}. 

The polytope $\Delta$ contains the Newton polytope of $F$: indeed, by taking a basis $\mbi{e}_1=(-5,1,1,0)$,\, $\mbi{e}_2=(1,1,-1,0)$,\, $\mbi{e}_3=(-3,0,0,1)$ for $\mathbb{R}^3$, one can see that monomials $WX^4$,\, $X^3Y$,\, $Y^2Z$,\, $Z^3$ are respectively corresponding to vertices 
\[
(1,2,-1),\, (1,1,-1),\, (0,-1,0),\, (-1,0,2). 
\]

The polytope $\Delta'$ contains the Newton polytope of $F'$: indeed, by taking a standard basis $\mbi{e}_1'=(-3,1,0,0)$,\, $\mbi{e}_2'=(-3,0,1,0)$,\, $\mbi{e}_3'=(-2,0,0,1)$ for $\mathbb{R}^3$, one can see that monomials $W'^9$,\, $X'Z'^3$,\, $X'^2Y'$,\, $Y'^3$ are respectively corresponding to vertices 
\[
(-1,-1,-1),\, (0,-1,2),\, (1,0,-1),\, (-1,2,-1). 
\]

The dual polytope $\Delta'^*$ of $\Delta'$ is a convex hull of vertices 
\[
(0,0,1),\, (-1,0,0),\, (-1,1,0),\, (0,1,0),\, (1,0,0),\,  (0,-1,-1)
\]
that is mapped to isomorphically from $\Delta$ by a transformation of $\mathbb{R}^3$ by the matrix 
\[
M=\left(
\begin{array}{ccc}
 2 & 2 & 1 \\
 -1 & 0 & 0 \\
 1 & 1 & 1 \\
\end{array}
\right)
\]
that is, $(x, y, z)M=(x', y', z')$ for $(x,y,z)\in\Delta$ and $(x',y',z')\in\Delta'$. 

Therefore, $\Delta$ and $\Delta'$ are reflexive and the pair is polytope dual. 

By the formula (\ref{SurjectionCriterion}), one gets $\rk L_0(\Delta)=\rk L_0(\Delta^*)=0$ because for all edges in $\Delta$ satisfy $l^*(\Gamma)l^*(\Gamma^*)=0$. 
In fact, at least either $\Gamma$ or $\Gamma^*$ has no lattice points in its interior. 

By the formula (\ref{Rho}), one can compute that 
\[
\rho(\Delta) = 20-3 = 17, \quad \rho(\Delta^*) = 6-3 = 3
\]
thus one has 
\[
\rho(\Delta) +\rho(\Delta^*) = 20. 
\]
\\

{\bf $Z_{2,0}$ and $Q_{17}$ case}
The defining polynomials of singularities $B=Z_{2,0}$ and $B'=Q_{17}$ are $f=x^5z +xy^3 + z^2,\, f'=x^5y +y^3+xz^2$, respectively. 

Take a compactification of $f$ as $F= W^7Y+X^5Z +XY^3 + Z^2$ in the weighted projective space $\mathbb{P}(1,1,3,5)$. 
Note that $F$ is the same compactification as in \cite{EbelingPloog}. 

Take a compactification of $f'$ as $F'=W^7+X^5Y +WY^3+XZ^2$ in the weighted projective space $\mathbb{P}(1,1,2,3)$. 
Note that $F'$ is the same compactification as in \cite{EbelingPloog}. 

The polytope $\Delta$ contains the Newton polytope of $F$: indeed, by taking a basis $\mbi{e}_1=(-3,3,0,0)$,\, $\mbi{e}_2=(-8,0,1,1)$,\, $\mbi{e}_3=(-6,1,0,1)$ for $\mathbb{R}^3$, one can see that monomials $W^7Y$,\,$X^5Z$,\, $XY^3$,\, $Z^2$ are respectively corresponding to vertices 
\[
(0,0,-1),\, (1,-1,1),\, (1,2,-3),\, (-1,-1,2). 
\]

The polytope $\Delta'$ contains the Newton polytope of $F'$: indeed, by taking a standard basis $\mbi{e}_1'=(-1,1,0,0)$,\, $\mbi{e}_2'=(-2,0,1,0)$,\, $\mbi{e}_3'=(-3,0,0,1)$ for $\mathbb{R}^3$, one can see that monomials $W^7,\, X^5Y,\, WY^3,\, XZ^2$ are respectively corresponding to vertices 
\[
(-1,-1,-1),\, (4,0,-1),\, (-1,2,-1),\, (0,-1,1). 
\]

The dual polytope $\Delta'^*$ of $\Delta'$ is a convex hull of vertices 
\[
(-1,-3,-4),\, (0,-2,-3),\, (0,1,0),\, (1,0,0),\, (0,0,1),\,  (-1,-2,-3)
\]
that is mapped to isomorphically from $\Delta$ by a transformation of $\mathbb{R}^3$ by the matrix 
\[
M := \left(
\begin{array}{ccc}
 1 & 1 & 1 \\
 1 & 3 & 4 \\
 1 & 2 & 3 \\
\end{array}
\right)
\]
that is, $M(x, y, z)=(x', y', z')$ for $(x,y,z)\in\Delta$ and $(x',y',z')\in\Delta'$. 

Therefore, $\Delta$ and $\Delta'$ are reflexive and the pair is polytope dual. 

By the formula (\ref{SurjectionCriterion}), one gets $\rk L_0(\Delta)=\rk L_0(\Delta^*)=0$ because for all edges in $\Delta$ satisfy $l^*(\Gamma)l^*(\Gamma^*)=0$. 
In fact, at least either $\Gamma$ or $\Gamma^*$ has no lattice points in its interior. 

By the formula (\ref{Rho}), one can compute that 
\[
\rho(\Delta) = 18-3 = 15, \quad \rho(\Delta^*) = 8-3 = 5
\]
thus one has 
\[
\rho(\Delta) +\rho(\Delta^*) = 20. 
\]
\\

{\bf $W_{1,0}$ case}
The defining polynomials of singularities $B=B'=W_{1,0}$ are the same $f=f'=x^6 +y^2z + z^2$. 

Take a compactification of $f$ as $F=X^6+Y^2Z+Z^2+W^6Z$ in the weighted projective space $\mathbb{P}(1,2,3,6)$. 
Note that $F$ is a different compactification from the one in \cite{EbelingPloog}. 

Take a compactification of $f'$ as $F'=X'^6+Y'^2Z'+Z'^2+W'^{12}$ in the weighted projective space $\mathbb{P}(1,2,3,6)$. 
Note that $F'$ is the same compactification as in \cite{EbelingPloog}. 

The polytope $\Delta$ contains the Newton polytope of $F$: indeed, by taking a basis $\mbi{e}_1=(-5,1,1,0)$,\, $\mbi{e}_2=(1,1,-1,0)$,\, $\mbi{e}_3=(-6,0,0,1)$ for $\mathbb{R}^3$, one can see that monomials $X^6$,\, $Y^2Z$,\, $Z^2$,\, $W^6Z$ are respectively corresponding to vertices 
\[
(2,3,-1),\, (0,-1,0),\, (-1,0,1),\, (-1,0,0). 
\]

The polytope $\Delta'$ contains the Newton polytope of $F'$: indeed, by taking a standard basis $\mbi{e}_1'=(-2,1,0,0)$,\,  $\mbi{e}_2'=(-3,0,1,0)$,\, $\mbi{e}_3'=(-6,0,0,1)$ for $\mathbb{R}^3$, one can see that monomials $X'^6$,\, $Y'^2Z'$,\, $Z'^2$,\, $W'^{12}$ are respectively corresponding to vertices 
\[
(5,-1,-1),\, (-1,1,0),\, (-1,-1,1),\, (-1,-1,-1). 
\]

The dual polytope $\Delta'^*$ of $\Delta'$ is a convex hull of vertices 
\[
(0,1,0),\, (-1,-1,-3),\, (0,-1,-2),\, (1,0,0),\, (0,0,1)
\]
that is mapped to isomorphically from $\Delta$ by a transformation of $\mathbb{R}^3$ by the matrix 
\[
M := \left(
\begin{array}{ccc}
 1 & 1 & 3 \\
 0 & 0 & -1 \\
 1 & 2 & 3 \\
\end{array}
\right)
\]
that is, $M(x, y, z)=(x', y', z')$ for $(x,y,z)\in\Delta$ and $(x',y',z')\in\Delta'$. 

Therefore, $\Delta$ and $\Delta'$ are reflexive polytopes and the pair is polytope dual. 

By the formula (\ref{SurjectionCriterion}), one gets $\rk L_0(\Delta)=\rk L_0(\Delta^*)=0$ because for all edges in $\Delta$ satisfy $l^*(\Gamma)l^*(\Gamma^*)=0$. 
In fact, at least either $\Gamma$ or $\Gamma^*$ has no lattice points in its interior. 

By the formula (\ref{Rho}), one can compute that 
\[
\rho(\Delta) = 21-3 = 18, \quad \rho(\Delta^*) = 5-3 = 2
\]
thus one has 
\[
\rho(\Delta) +\rho(\Delta^*) = 20. \qquad\textnormal{\QED } 
\]

\hfill \textsc{Makiko Mase} \\ 
\hfill e-mail: mtmase@arion.ocn.ne.jp \\
\hfill {\small{Department of Mathematics, Rikkyo University}} \\
\hfill {\small{171-8501 3-34-1 NishiIkebukuro, Toshimaku, Tokyo, Japan. }}\\
\end{document}